# ISOMORPHIC PASTINGS AND THE TWO POSSIBLE STRUCTURES FOR A PAIR OF GRAPHS HAVING THE SAME DECK

S.Ramachandran
Mathematics Dept.,
Vivekananda College,
Agasteeswaram-629701, T.N., INDIA
dr.s.ramachandran@gmail.com

**Abstract**

The unlabeled subgraph obtained by deleting a vertex from a graph G is a card of G and the collection of cards of G is the deck of G. A graph having the same deck as G is a hypomorph of G. A graph is called reconstructible if it is isomorphic to all its hypomorphs. Reconstruction Conjecture claims that all graphs are reconstructible and it is open. A representation of a hypomorph of G in terms of two of its cards, called pasting, is introduced and studied. Isomorphic pastings of two cards are defined. In the case of a digraph, a card with which the degree triple of the deleted vertex is also given is called a degree associated card or dacard. Dadeck, da-reconstructible digraphs, dapastings etc. are defined analogously. DARC claims that all digraphs are da-reconstructible and it is also open. *Results*: Two hypomorphs G and H of a graph are isomorphic if and only if a pair of cards in their common deck is pasted isomorphically in both G and H. Either every pair of cards in their common deck is pasted isomorphically in both G and H, or no pair of cards is pasted isomorphically in both G and H. Results similar to the above hold for dapastings in da-hypomorphs of a digraph. Two new graph parameters, the neighborhood degree quintuple of a vertex and a new family of digraphs are da-reconstructible. New approaches for proving the reconstruction conjecture and DARC by the method of contradiction arise.

## 1. Introduction

We use the terminology in Harary ([5]) for graphs/digraphs and in [3, 7] for reconstruction. In a digraph D, a vertex $w$ is called an *out-neighbor*, *in-neighbor* or *strong neighbor* of a vertex $v$ according as $vw$ is an unpaired arc, $wv$ is an *unpaired arc* or $vw$ and $wv$ are both arcs. When $vw$ and $wv$ are both arcs, they together are called a *biarc* or *a symmetric pair of arcs*; when none of them is an arc, $v$ and $w$ are called non-adjacent. When $a$, $b$ and $c$ denote the numbers of out-neighbors, in-neighbors and strong neighbors of $v$, the ordered triple $(a, b, c)$ is called the *degree triple* of $v$, denoted as $dt(v)$.

**Definition 1.1.** For a vertex $v$ of a graph (digraph) G, the unlabeled subgraph (subdigraph) $G-v$ is called a *card* of G and the multiset $\{G-v \mid v \in V(G)\}$ is called the *deck* of G ([4]) and is denoted by *Deck*(G). When $G-v$ is a card of a graph (digraph) G, the ordered pair $(G-v, d_G(v))$ $((G-v, dt_G(v)))$ is called a *degree associated card* or *dacard* of G and the collection of dacards of G is called the *dadeck* of G ([8, 7]) denoted as *Dadeck*(G). Graphs/digraphs having the same deck (dadeck) are called *hypomorphs* (*da-hypomorphs*) of each other. A graph/digraph which is isomorphic to all its hypomorphs (da-hypomorphs) is called *reconstructible* ([4]) (*da-reconstructible* ([8])).

**Observation 1.2.** In the case of a graph G with $|G| \geq 3$, $d(u)$ can be determined from card G – u and *Deck*(G) ([3]) and hence *Deck*(G) gives *Dadeck*(G).

**Reconstruction Conjecture** (**RC**) ([14, 6]). If G and H are two hypomorphic graphs each on at least three vertices, then G ≅ H. (That is, all graphs on ≥ 3 vertices are reconstructible.)

Harary in 1964 *extended* RC to digraphs ([4]) as follows.

**Digraph Reconstruction Conjecture (DRC).** All digraphs on ≥ 5 vertices are reconstructible.

Stockmeyer in 1977 exhibited infinite numbers of pairs of digraphs ([11]) that disobey DRC and so a weaker form, which is also an extension of RC was proposed and both are still open. Work being done on RC and related problems are surveyed in [3, 7].

**Degree Associated Reconstruction Conjecture** (**DARC**) ([8]). All digraphs on ≥ 2 vertices are da-reconstructible.

It has been proved in [1] that each graph/digraph in some infinite families is uniquely determined by two of its dacards. Reconstruction from *neighborhood degrees associated cards* (G–v, ND(v)) is studied in [9]. Here we give a *representation* of a hypomorph of an arbitrary graph (digraph) in terms of *two of its cards* (*dacards*) *for the purpose of studying isomorphism between two hypomorphs*. This key concept is called *a pasting of two cards* (*a dapasting of two dacards*).

**Definition 1.3** ([10]). Given a digraph (respectively, graph) G and two distinct vertices $v_1$ and $v_2$, let A ≅ G – $v_1$ and B ≅ G – $v_2$. Suppose P is a digraph (respectively, graph) with two distinguished non-adjacent vertices $u$ and $v$, both labeled with $e$, such that

(a). P – u ≅ A and P – v ≅ B, and

(b). there exists a digraph H ∈ {P, P + uv, P + vu, P + uv + vu} (respectively, graph H ∈ {P, P + uv}) such that H is hypomorphic to G.

Then (P, Deck(G)) is *a pasting of A and B as members of Deck*(G) and H is *a completion* of this pasting. Vertices of P labeled with $e$ are *the external vertices of the cards* in the pasting.

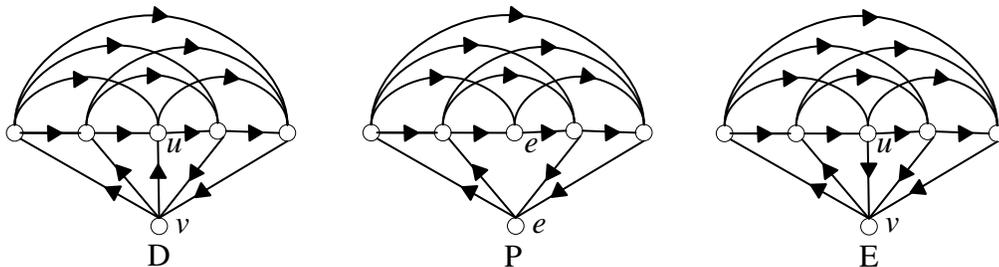

**Figure 1.** Two hypomorphic digraphs D and E with D ≇ E and a pasting P of D – u and D – v.

**Example 1.4.** The pair (D, E) of hypomorphic but nonisomorphic digraphs in Figure 1 was discovered by Beineke and Parker ([2]) and is given in [5, p. 208] also. The digraph P is a pasting of the cards D – u and D – v as members of *Deck*(D) and D is a completion of P. Again P is a pasting of the cards E – u and E – v of E as members of *Deck*(E) and E is also a completion of P.

A definition similar to Definition 1.3 holds for *a dapasting of two dacards of a digraph as members of Dadeck(G)*. We use these definitions together with a notion of isomorphism between pastings to refine our knowledge of reconstruction of graphs and da-reconstruction of digraphs. Properties of non-reconstructible graph/digraph pairs in general have been studied in [13] and known pairs of non-reconstructible digraphs in [12, 13, 10] and some others. It is known from 1957 itself ([6]) that the claim of RC is in terms of a pair of graphs having the same deck. The structure of such a pair of graphs in terms of cards from their common deck is found only now (in Theorem 5.11 below). Also in Theorem 5.8, we prove that two graphs (digraphs) having the same deck (dadeck) are isomorphic if and only if a pair of cards (dacards) in their common deck (dadeck) is pasted (dapasted) isomorphically in both of them. Some relations connecting the parameters of a dapasting and the da-reconstructible parameters of a graph/digraph are derived. Such results may even facilitate proving that RC (DARC) cannot have a counterexample pair.

## 2. Definitions and notations

The number of vertices in a graph G is denoted by $|G|$. Two vertices $u$ and $v$ of G are called *similar* if there is an automorphism of G taking $u$ to $v$. Similarity is an equivalence relation on V(G) and partitions it into similarity classes called *orbits* of G. The number of subgraphs of G isomorphic to F is denoted by $S(F, G)$, the number among them containing $v$ by $S(F, G, v)$ and the number among them containing $v$ in the orbit $\theta$ of the subgraph F by $S(F, G, v, \theta)$. The corresponding numbers for induced subgraphs are denoted by $s(F, G)$, $s(F, G, v)$ and $s(F, G, v, \theta)$ respectively. The above definitions/notations *apply for digraphs also*. When $dt(v) = (a,b,c)$, the numbers $a$, $b$, $c$, $a+c$ and $b+c$ are called the *first degree*, *second degree*, *third degree*, *outdegree* ($od(v)$) and the *indegree* ($id(v)$) of $v$ respectively ([5]). The ordered pairs ($od(v)$, $id(v)$) and ($a+b$, $c$) are called the *degree pair* ([5]) and the *split degree of v* (written as $dep(v)$ and $spd(v)$) respectively). When G is a graph, the collection of degrees of the neighbors of $v$ is called the neighborhood degree sequence ($ND_G(v)$ or $ND(v)$) of $v$.

**Definition 2.1.** For each vertex $v$ of a digraph D, the ordered quintuple ($csdon(v)$, $cidon(v)$, $cfdin(v)$, $codin(v)$, $ctdsn(v)$) is called the *neighborhood degree quintuple of v* (denoted as $ndq(v)$) where $csdon(v)$ denotes the collection of second degrees of the out-neighbors of $v$ ([9]),

$cidon(v)$ denotes the collection of in-degrees of the out-neighbors of $v$,

$cfdin(v)$ denotes the collection of first degrees of the in-neighbors of $v$ ([9]),

$codin(v)$ denotes the collection of out-degrees of the in-neighbors of $v$ and

*ctdsn(v)*) denotes the collection of third degrees of the strong neighbors of *v* ([9]).

**Definition 2.2.** A parameter or a function defined on graphs/digraphs G of a class $\mathbb{C}$ of graphs/digraphs is called a *da-reconstructible parameter of a graph/digraph* of class $\mathbb{C}$ if it takes the same value on all da-hypomorphs of G whenever G ∈ $\mathbb{C}$ [7, 8, 9]. A parameter P of a vertex *v* of a digraph G in a class $\mathbb{C}$ of digraphs is called a *da-reconstructible parameter of a vertex* (or *of a dacard*) for *class* $\mathbb{C}$ if P(*v*) = P(*v**) whenever G ∈ $\mathbb{C}$, H is a da-hypomorph of G, G–*v* ≅ H–*v** and $dt_G(v) = dt_H(v^*)$. *Reconstructible parameters* are defined analogously ([4]).

**Notation 2.3.** The number of induced paths of length *k* in graph G starting at *v* is denoted as *ip*(*k*, G, *v*). The number of induced *u*–*v* paths of length *k* in G is denoted as *ip*(*k*, G, *u*, *v*).

**Notation 2.4.** If pastings of two cards of G considered are always pastings as members of *Deck*(G) in a discussion as *in this article*, then "a pasting of A and B as members of *Deck*(G)" is called "a *pasting* of A and B" in short form and written as P.

## 3. Some reconstructible parameters whose effect on pastings is easily expressible

Most of the parameters known to be reconstructible are listed in [3, 7]. Here we reconstruct some new parameters and restate a few in forms in which they are used later.

**Lemma 3.1. (a). (Kelly's lemma** [6, 3, 10]**).** For all graphs/digraphs F with |F| < |G|, S(F, G) and *s*(F, G) are reconstructible and S(F, G, *v*) and *s*(F, G, *v*) are determinable from G–*v* and *Deck*(G).

**(b).** The number of unpaired arcs and the number of biarcs in a digraph D are reconstructible.

**(c).** For each vertex *v* of a digraph D, $spd_D(v)$ is determinable from D–*v* and *Deck*(D). If $dep_D(v)$ is given, then $dt_D(v)$ and (D–*v*, $dt_D(v)$) are determinable from D–*v* and *Deck*(D). Thus *Dadeck*(D) is determinable from the collection {(D–*v*, $dep_D(v)$) | *v* ∈ V(D)}.

**Proof.** Parts (b) and (c) can be easily proved from Part (a) and is omitted. □

**Lemma 3.2.** When G is a graph and *v* ∈ V(G), ND(*v*) is determinable ([3]) from G–*v* and *Deck*(G).

**Lemma 3.3. (a)** When F is a star with *m* vertices satisfying 3 ≤ *m* < |G| and θ is an orbit of F, S(F, G, *v*, θ) is a reconstructible parameter of a vertex *v* of a graph G. In addition, if G is a *triangle-free graph*, then *s*(F, G, *v*, θ) is also a reconstructible parameter of a vertex *v*.

**(b)** For a graph G with |G| ≥ 3, *ip*(2, G, *v*) is a reconstructible parameter of a vertex *v*.

**Proof.** (a) F has exactly two orbits. One of them, say $θ_1$ has only the center of F in it and the other, say $θ_2$ contains all the end vertices of F. By Kelly's lemma S(F, G, *v*) is known. Since m ≥ 3, we have S(F, G, *v*, $θ_1$) = C($deg_G(v)$, *m*–1) (where C(*n*, *r*) denotes the number of selections of *r* objects from *n* distinct objects) and so is determinable from *Deck*(G) and the card G–*v*. Hence S(F, G, *v*, $θ_2$) = S(F, G, *v*) – S(F, G, *v*, $θ_1$) is also a reconstructible parameter of a vertex of G. (1)

When G is *triangle-free*, every star on at least three vertices is an induced star and so $s(F, G, v, \theta) = S(F, G, v, \theta)$ for all orbits $\theta$ of F and hence is a reconstructible parameter of a vertex of G.

(b) The required number $ip(2, G, v) = s(K_{1,2}, G, v, \theta_2) = S(K_{1,2}, G, v, \theta_2) - 2S(K_3, G, v)$ since each $K_3$ containing $v$ contributes two $K_{1,2}$'s to $S(K_{1,2}, G, v, \theta_2)$, which are not to be counted in $s(K_{1,2}, G, v, \theta_2)$. Since $S(K_{1,2}, G, v, \theta_2)$ and $S(K_3, G, v)$ are determinable from *Deck*(G) and G–v by Part (a) above and Kelly's Lemma respectively, the claim follows. □

**Definition 3.4** ([9]). When $w$ is a vertex of digraph D, the ordered triple $(C_1, C_2, C_3)$ where $C_1$, $C_2$ and $C_3$ are the collections of degree triples of the out-neighbors, in-neighbors and strong neighbors of $w$ in digraph D is called the *neighborhood degree triple* (NDT) of $w$ in D and is denoted as $NDT_D(w)$.

**Observation 3.5.** (a) The multiset of degree triples of the vertices of a digraph G is da-reconstructible. (b) The collection $\{(dt(v), \text{collection of NDT's of the vertices of } D-v)| v \in V(D)\}$ is da-reconstructible.

For graphs G, the neighborhood degree sequence of a vertex $v$, $(ND(v))$, is an important reconstructible parameter ([3]). For da-reconstruction of digraphs, we prove the following.

**Lemma 3.6.** For all digraphs D, $ndq(v)$ is a da-reconstructible parameter of a vertex $v$.

**Proof.** We show that each of the parameters $csdon(v)$, $cidon(v)$, $cfdin(v)$, $codin(v)$, $ctdsn(v)$ are determinable from the dadeck of D and the dacard D–v.

Determination of $cidon(v)$ and $codin(v)$: The degree pair sequence of D is known from the dadeck of D. Also for each dacard D–v, $dt_D(v)$ is known. "The number of out-neighbors of $v$ whose indegree in D is $k$" is the difference between "the number of vertices which have indegree less than $k$ in the digraph D–v" and the "number of vertices other than $v$ whose indegree in D is less than $k$". This is reconstructible as the last two numbers are da-reconstructible. Thus the number of out-neighbors of $v$ whose indegree in G is $k$ can be determined for each $k$, $1 \leq k \leq n-1$ and hence $cidon(v)$ can be formed from them. By the same method, $codin(v)$ can be found from dacard D–v and *Dadeck*(D). The three parameters $csdon(v)$, $cfdin(v)$ and $ctdsn(v)$ are determined from dacard D–v and *Dadeck*(D) using similar arguments in [9, Theorem 2.2]. □

**Observation 3.7.** The six *collections* of ordered pairs $\{(dt(v), \text{second degree of } w) \mid vw \text{ is an unpaired arc of D}\}$, $\{(\text{first degree of } v, dt(w)) \mid vw \text{ is an unpaired arc of D}\}$, $\{(\text{first degree of } v, \text{second degree of } w) \mid vw \text{ is an unpaired arc of D}\}$, $\{(dt(v), id(w) \mid vw \text{ is an unpaired arc of D}\}$, $\{(od(v), dt(w)) \mid vw \text{ is an unpaired arc of D}\}$ and $\{(dt(v), ndq(v))| v \in V(D)\}$ are da-reconstructible parameters of a digraph. (This is used in the proof for Observation 7.5 below.)

## 4. Pastings and dapastings

We first define a dapasting of two dacards.

**Definition 4.1.** Let G be a digraph and let (A, $dt(A)$) and (B, $dt(B)$) be two dacards of G obtained by deleting distinct vertices of G (i.e., there are distinct vertices $w_1$ and $w_2$ of G such that $A \cong G - w_1$, $B \cong G - w_2$, $dt(A) = dt_G(w_1)$ and $dt(B) = dt_G(w_2)$). Suppose P is a digraph with two distinct vertices $u$ and $v$ labeled with $(e, dt(A))$ and $(e, dt(B))$ respectively, such that

(a). $P - u \cong A$ and $P - v \cong B$, and

(b). there exists a digraph H da-hypomorphic to G such that H is isomorphic to Y for some Y ∈ {P + uv, P + vu, P + uv + vu, P}, $dt_Y(u) = dt(A)$ and $dt_Y(v) = dt(B)$.

Then (P, *Dadeck*(G)) is *a dapasting of the dacards* (A, $dt(A)$) *and* (B, $dt(B)$) *as members of Dadeck*(G) and H is *a completion* of the dapasting (P, *Dadeck*(G)). The vertices $(e, dt(A))$ and $(e, dt(B))$ in P are called *the external vertices of the dacards* (A, $dt(A)$) and (B, $dt(B)$) respectively in the dapasting P.

**Remark.** The above definition gives the following. When $u$ and $v$ are vertices of a digraph G and P denotes the digraph obtained from G by deleting all arcs between $u$ and $v$ and labeling $u$ and $v$ with $(e, dt_G(u))$ and $(e, dt_G(v))$ respectively, the ordered pair (P, *Dadeck*(G)) will be a dapasting of the dacards $(G-u, dt_G(u))$ and $(G-v, dt_G(v))$ of G as members of *Dadeck*(G).

**Definition 4.2.** A dapasting (P, *Dadeck*(G)) of dacards (A, $dt(A)$) and (B, $dt(B)$) of G is said to be a dapasting of (A, $dt(A)$) and (B, $dt(B)$) *in a da-hypomorph* J *of* G if J is a completion of P. A pasting P of cards A and B *in a hypomorph* J *of* G is defined similarly. (See Example 4.11 and Lemma 4.12 also.)

As digraphs D and E in Figure 1 above are hypomorphic and E is a completion of P, P is a pasting of the cards D – $u$ and D – $v$ as members of *Deck*(D) in E also. Similarly P is a pasting of the cards E – $u$ and E – $v$ as members of *Deck*(E) in D.

**Observation 4.3. (a).** While studying a pasting (dapasting) P of two cards (dacards) of G as members of *Deck*(G) (*Dadeck*(G)), in addition to P we know *Deck*(G) (*Dadeck*(G)) also and hence we can derive many relations (Observation 4.14, for example) connecting the parameters of P with the reconstructible parameters of G and of the two cards that are pasted (dapasted) in P.

**(b).** Every pasting (dapasting) has at least one completion and if H is a completion of a pasting (dapasting) P, then all graphs/digraphs isomorphic to H are also completions of P by Definitions 1.3(b) and 4.1(b).

**(c).** As digraphs without unpaired arcs can be taken as graphs, Definition 4.1 with corresponding simplifications can be taken to define dapasting of a pair of dacards of a graph also. In the case of a graph G, a pasting of a pair of cards of G as members of *Deck*(G) gives a dapasting of the corresponding dacards as members of *Dadeck*(G) and vice versa (by Observation 1.2).

*In this paper a pasting (dapasting) of two cards (dacards) of G is considered to be as members of Deck(G) (Dadeck(G)) unless stated otherwise. The study is done for the dapasting of dacards of a digraph as it covers dapastings of dacards of a graph (which are in fact pastings of cards of a graph) as a particular case* (Observation 4.3(c)). *We point out when the corresponding result for pastings of cards of a digraph is not true.*

The following is a direct consequence of Definitions 1.3, 4.1 and Lemma 3.1(b).

**Observation 4.4.** When P is a dapasting of two dacards (A, $dt$(A)) and (B, $dt$(B)) of a digraph G as members of *Dadeck*(G) with $u$ ($v$) as the external vertex of A (B), the following hold.

**(a).** *At least one of the digraphs* H "obtained by joining $u$ and $v$ in P by means of arc(s) (if needed) such that the degree triples of $u$ and $v$ in H become $dt$(A) and $dt$(B) respectively (that is, attain the values given in the second co-ordinates of their labels in P)" is a da-hypomorph of G with dacard (A, $dt$(A)) $\cong$ (H−$u$, $dt_H(u)$) and dacard (B, $dt$(B)) $\cong$ (H−$v$, $dt_H(v)$). Also, *each completion* H of P can be obtained from P by the above process.

**(b).** A *vertex proper subdigraph* of P must be a subdigraph of at least one of the cards of G ([10, Theorem 3.4]).

**(c).** When P is a pasting of two cards $A_1$ and $A_2$ of a graph (digraph) G with $u_1$ and $u_2$ as the external vertices of the cards such that P − $u_i$ $\cong$ $A_i$, $i$ = 1, 2, $u_i$ can be labeled with ($e$, $d_i$) where $d_i$ is the degree (split degree) of a vertex $u_i^*$ in a completion H of P such that H − $u_i^*$ $\cong$ $A_i$.

**Observation 4.5.** A dapasting is a union of a chosen pair of dacards of D regulated by the dadeck as whose members they are considered. So, whether a digraph/graph is a dapasting of two given dacards or not depends on the dadeck with respect to which we seek the answer. Hence the expression *as members of Dadeck*(G) in Definition 4.1 is an important part as illustrated below.

**Example 4.6.** In Figure 2, let A $\cong$ G–$u$ $\cong$ R–$a$, and B $\cong$ G–$v$ $\cong$ R–$b$. "G with $u$ and $v$ relabeled as ($e$, 1)" is a dapasting of dacards (A, 1) and (B, 1) as members of *Dadeck*(G), but it is not a dapasting of them as members of *Dadeck*(R); because "G cannot be a subgraph of a hypomorph of R as required by Definition 4.1 (b). (The graph F, which is a subgraph of G with |F| < |R| is not a subgraph of R and hence by Kelly's lemma cannot be a subgraph of a hypomorph of R.) Similarly, R with $a$ and $b$ relabeled as ($e$, 1) is a dapasting of (A, 1) and (B, 1) as members of *Dadeck*(R), but it is not a dapasting of them as members of *Dadeck*(G); because, "R cannot be a subgraph of a da-hypomorph of G" (as the subgraph T of R with |T| < |G| is not a subgraph of G).

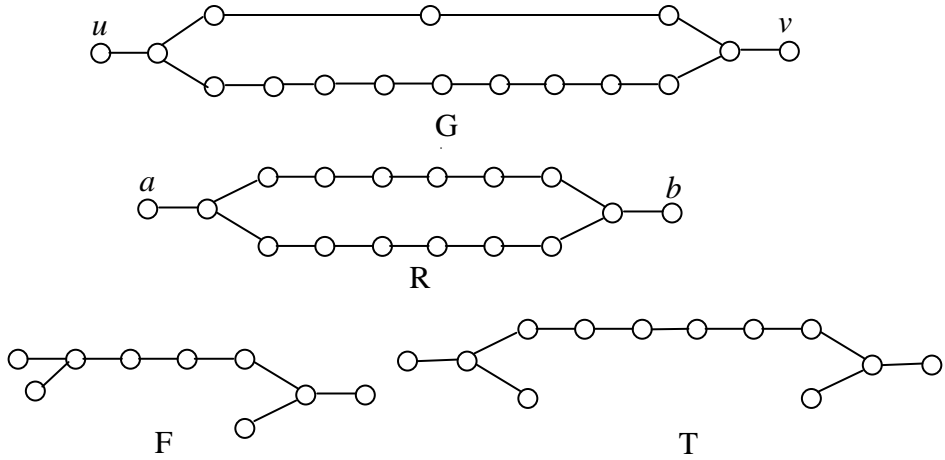

**Figure 2.** G–u ≅ R–a, G–v ≅ R–b and F (T) is a vertex proper subgraph of G (R), but not of R (G).

**Lemma 4.7.** Let G be a digraph and let $A_i = G - u_i$, $i = 1, 2$, where $u_1$ and $u_2$ are distinct vertices of G. (a). *Every da-hypomorph H of G can be obtained as a completion of a dapasting of the dacards $(A_1, dt_G(u_1))$ and $(A_2, dt_G(u_2))$.*

(b). *Every hypomorph H of G can be obtained as a completion of a pasting of $A_1$ and $A_2$.*

**Proof.** (a): Since H is a da-hypomorph of G, there exist distinct vertices $w_1$ and $w_2$ in H such that $H - w_i \cong A_i$ and $dt_H(w_i) = dt_G(u_i)$ for $i = 1, 2$. Now the dacards $(H - w_i, dt_H(w_i))$ and $(A_i, dt_G(u_i))$ are identical for $i = 1, 2$. Hence by Definition 4.1, the digraph P obtained from H by deleting all arcs between $w_1$ and $w_2$ and labeling $w_i$ with $(e, dt_G(u_i))$ for $i = 1, 2$ is a dapasting of the dacards $(A_1, dt_G(u_1))$ and $(A_2, dt_G(u_2))$ of G as members of *Dadeck*(G) and H is a completion of it.

Claim (b) can be proved similarly. □

**Definition 4.8** ([10])**.** Two dapastings (pastings) P1 and P2 of a pair of *dacards* (*cards*) of a digraph/graph G as members of *Dadeck*(G) (*Deck*(G)) are called *isomorphic* if there is an isomorphism from P1 to P2 *preserving the labels*.

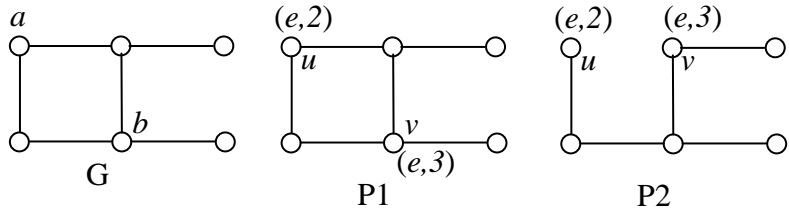

**Figure 3:** A graph G and two nonisomorphic dapastings P1 and P2 of a pair of dacards isomorphic to G–a and G–b in G.

**Example 4.9.** There may be more than one (nonisomorphic) dapastings of a given pair of dacards of G. For the graph G in Figure 3, P1 and P2 are two nonisomorphic dapastings of the dacards (G

– a, 2) and (G – b, 3) as members of *Dadeck*(G). Both of them occur as dapastings in G with vertices labeled (*e*, 2) and (*e*, 3) as the external vertices of the cards.

**Lemma 4.10.** Let P be a pasting (da-pasting) of two cards (dacards) of a digraph G in a hypomorph (da-hypomorph) of G. *The number and the type of arcs* with which the external vertices of the cards in P are to be joined to get a completion of P are unique and are determined by P and G. The corresponding results hold for the pastings (dapastings) of two cards (dacards) of a graph also.

**Proof.** The number of unpaired arcs and the number of biarcs in a digraph G are reconstructible parameters of a digraph (Lemma 3.1(b)). These respectively give the number of unpaired arcs and the number of biarcs in the completion of P since the completion is always a hypomorph of G (by Definition 1.3(b) and Definition 4.1(b)). These on comparison with the corresponding values for P give the exact number and the type of arc(s) with which the external vertices of A and B in P are to be joined in order to get a completion of P. □

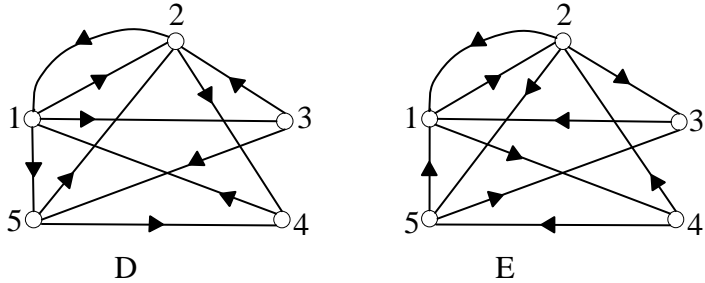

**Figure 4.** Two nonisomorphic hypomorphs with a pasting in one isomorphic to a subdigraph, (but not a pasting) in the other.

**Example 4.11.** The digraphs D and E in Figure 4 form a *non-reconstructible pair* ([10, Fig. 1]) and $D - j \cong E - j$ for $j = 4, 5$. Also D – arc (5, 4) with vertices 4 and 5 relabeled "*e*" is *a pasting P of cards* D – 4 and D – 5 in D. Also P is isomorphic to the subdigraph E – arc (4, 5) of E under an isomorphism β taking vertices 1, 2, 3, **4** and **5** of D to vertices 2, 1, 5, **4** and **3** of E respectively. But the subdigraph β(P) ($\cong$ P) with β(4) (= 4) and β(5) (= 3) labeled "*e*" *is not a pasting of the cards* D – 4 and D – 5 *as members* of *Deck*(D) *in* E since E is *not isomorphic to* any hypomorph of D obtained by joining the vertices 4 and 3 of β(P) (which correspond to the external vertices of P) by means of arc(s) and so E is not a completion of β(P). (We have β(P) + arc(3, 4) $\cong$ D $\not\cong$ E and β(P) + arc(4, 3) $\not\cong$ E (as degree triples of vertices in them differ).

If the above situation happens in the case of two da-hypomorphs of a digraph (graph), they must be nonisomorphic and disprove DARC (RC) as the following lemma shows.

**Lemma 4.12.** If G and H are da-hypomorphic graphs/digraphs and there exists a dapasting P of two dacards of G as members of *Dadeck*(G) such that P is a subgraph/subdigraph of H, but not a dapasting in H, then G is not da-reconstructible.

**Proof.** P is a dapasting of two dacards of G as members of *Dadeck*(G). Hence the completion CP of P is a da-hypomorph of G. As P is not a dapasting in H, CP $\not\equiv$ H (by Definition 4.2). But CP and H are both da-hypomorphs of G and hence G is not da-reconstructible. □

The following lemma is used in the proof of the two results that succeed it.

**Lemma 4.13.** Let $u$ and $v$ be vertices of a digraph G and let P be the digraph obtained from G by deleting all arcs (if any) joining $u$ and $v$. If the arcs joining $u$ and $v$ in G are known, then $dt_P(u)$, $dt_P(v)$, $ndq_P(u)$ and $ndq_P(v)$ can be determined in terms of $dt_G(u)$, $dt_G(v)$, $ndq_G(u)$ and $ndq_G(v)$.

**Proof.** Obviously, $dep(w) = (od(w), id(w))$ can be determined from $dt(w)$ for any vertex $w$. There are four ways in which $u$ and $v$ can be joined in G.

**Case 1.** *Let vu be the only arc joining u and v in G*.

All out-neighbors of $u$ in G occur as out-neighbors of $u$ in P with the same degree triple as in G. Among the in-neighbors of $u$ in G, only $v$ ceases to be an in-neighbor of $u$ in P and others occur as in-neighbors of $u$ in P with the same degree triple as in G. Hence only $dt(u)$ and the components corresponding to $cfdin(u)$ and $codin(u)$ in $ndq(u)$ will change in P. Thus "$dt_P(u) = dt(A)–(0,1,0)$", "$cfdin_P(u) = cfdin_G(u) - \{$first entry of $(dt(B))\}$", "$codin_P(u) = codin_G(u) - \{$first entry of $(dep(B))\}$", "$(csdon(u), cidon(u), ctdsn(u))$ in P and in G are equal".

The values of these parameters for $v$ are given by "$dt_P(v) = dt(B) - (1,0,0)$", "$csdon_P(v) = csdon_G(v) - \{$second entry of $(dt(A))\}$", "$cidon_P(v) = cidon_G(v) - \{$second entry of $(dep(A))\}$" and "$(cfdin(v), codin(v), ctdsn(v))$ in P and in G are equal".

*In the other three cases also*, the values of these parameters can be computed similarly. □

The following two results are examples of the ways in which the parameters of a da-pasting of two dacards of a graph/digraph G are constrained by the da-reconstructible parameters of G.

**Observation 4.14.** Let G be a digraph and let (A, $dt(A)$) and (B, $dt(B)$) be dacards of G obtained by deleting distinct vertices of G. If P is *a dapasting* of A and B as members of *Dadeck*(G) with $u$ and $v$ as the external vertices of the dacards A and B respectively, then the following hold.

(a). There exist isomorphisms $\beta_A$: A $\to$ P–$u$ and $\beta_B$: B $\to$ P–$v$ and vertices $u^* \in V(B)$ and $v^* \in V(A)$ such that $\beta_B^{-1}\beta_A$ (function $\beta_A$ acts first) is an isomorphism from A – $v^*$ to B – $u^*$.

(b). The multiset of degree triples of the vertices of digraph G after subtraction of the multiset $\{dt(A), dt(B)\}$ from it is same as the multiset of degree triples of the vertices of the digraph P other than $u$ and $v$.

(c). The values of $dt_P(u)$, $dt_P(v)$, $ndq_P(u)$ and $ndq_P(v)$ can be expressed in terms of $dt_G(u)$, $dt_G(v)$, $ndq_G(u)$ and $ndq_G(v)$.

**Proof.** Part (a) follows directly from Definition 4.1.

*Parts* (b) *and* (c): Let $u$ and $v$ be distinct vertices of G such that $(G-u, dt_G(u)) \cong (A, dt(A))$ and $(G-v, dt_G(v)) \cong (B, dt(B))$. By Definition 4.1(b), there exists a digraph H, da-hypomorphic to G such that either (i) $H \cong P$ or (ii) $H \cong P + uv + vu$ or (iii) $H \cong P + uv$ or (iv) $H \cong P + vu$, satisfying $(A, dt(A)) \cong (H-u, dt_H(u))$ and $(B, dt(B)) \cong (H-v, dt_H(v))$. As G and H are da-hypomorphs, by Observation 3.5 (a) and Lemma 3.6, G and H have the same collection of degree triples for their vertices, $dt_H(u) = dt_G(u)$, $dt_H(v) = dt_G(v)$, $ndq_H(u) = ndq_G(u)$ and $ndq_H(v) = ndq_G(v)$. (1)

Obviously, P is the digraph obtained from H by deleting all arcs (if any) joining $u$ and $v$, and hence Part (b) now follows from (1). Moreover, by Lemma 4.13, in each of the four situations given above, the values of $dt_P(u)$, $dt_P(v)$, $ndq_P(u)$ and $ndq_P(v)$ can be expressed in terms of $dt_H(u)$, $dt_H(v)$, $ndq_H(u)$ and $ndq_H(v)$ and hence, Part (c) now holds by (1). □

In the case of pastings of cards of a digraph (graph), instead of degree triple, only the split degree (degree) of the missing vertex is proved to be determinable for each card (Lemma 3.1(c)) using the deck and so we have the following.

**Lemma 4.15.** Observation 4.14 gives some properties of a *pasting of a pair of cards* of a digraph (graph) when the degree triple of the missing vertex associated with a dacard is replaced by the split degree (degree) of the missing vertex and conditions from (a) to (c) modified accordingly. When G is a graph, P satisfies condition (d) also.

(d). (i) If $|E(G)| = |E(P)|$, then $ip(2, G, u) = ip(2, P, u)$, $ip(2, G, v) = ip(2, P, v)$, $S(K_3, G, u)$
$= S(K_3, P, u)$ and $S(K_3, G, v) = S(K_3, P, v)$. (In all the above, L.H.S. is reconstructible.)

(ii) If $|E(G)| = |E(P)| + 1$, then $ip(2, G, u) = ip(2, P, u) - 2ip(2, P, u, v) + d(B) - 1$,

$$ip(2, G, v) = ip(2, P, v) - 2ip(2, P, v, u) + d(A) - 1,$$

$$S(K_3, G, u) = S(K_3, P, u) + ip(2, P, u, v) \text{ and}$$

$$S(K_3, G, v) = S(K_3, P, v) + ip(2, P, v, u).$$

**Proof.** (i) Now P is a hypomorph of G with $A \cong P - u$ and $B \cong P - v$ and the claims are obvious.

(ii) Now $P + uv$ is a hypomorph of G with $A \cong (P + uv) - u$ and $B \cong (P + uv) - v$. Hence, $ip(2, G, u)$ = (number of those counted in $ip(2, G, u)$ which do not contain $v$) + (number of those counted in $ip(2, G, u)$ which contain $v$) = $(ip(2, P, u) - ip(2, P, u, v)) + (d(B) - 1 - ip(2, P, u, v))$. Others can be proved similarly. (L.H.S. is found from deck and R.H.S. is observed in the pasting P). □

**5. Isomorphism between two hypomorphs and isomorphism between pastings in them**

The following can be easily derived from Definition 4.8.

**Lemma 5.1.** If $P_j$ is a dapasting of the dacards ($A_j$, $dt(A_j)$) and ($B_j$, $dt(B_j)$) of a digraph G as members of *Dadeck*(G) for $j = 1, 2$ and $P_1$ and $P_2$ are isomorphic, then the two collections {($A_1$, $dt(A_1)$), ($B_1$, $dt(B_1)$)} and {($A_2$, $dt(A_2)$), ($B_2$, $dt(B_2)$)} of dacards of G are the same. The corresponding result for isomorphic pastings of graphs/digraphs also holds.

**Theorem 5.2.** Let G be a graph or a digraph and let A and B be dacards obtained by deleting distinct vertices of G. Let P1 and P2 be two dapastings of A and B as members of *Dadeck*(G). There is "an isomorphism between the completions of P1 and P2 *taking the set of external vertices of the cards in one dapasting to that in the other*" if and only if the two dapastings are isomorphic.

**Proof.** Let CP1 be a completion of P1 and CP2 be a completion of P2. *We prove the theorem for dapastings of digraphs and the proof for dapastings of graphs will follow.*

*"Only if" part*: By hypothesis there is an isomorphism φ from CP1 to CP2 taking the set of external vertices of P1 to that of P2. Since CP1 gives P1 (CP2 gives P2) when the arcs between the external vertices of the cards in it are deleted (by the definition of completion of a dapasting), φ is an isomorphism from P1 to P2 taking the set of external vertices of P1 to that of P2.

*"If" part*: Let β be an isomorphism from P1 to P2. Let $u$ ($w$) denotes the external vertex of A (B) in P1. Hence by Definition 4.8, {β($u$), β($w$)} is the set of external vertices of A and B in P2.   (1)

Since the completions of both P1 and P2 are da-hypomorphs of G, they have the same number of unpaired arcs and the same number of biarcs as G. Hence as P1 $\cong$ P2, the type and the number of arcs to be added (if any) between the external vertices of A and B while completing P1 are same as those of the arcs to be added between β($u$) and β($w$) while completing P2 and they are known from P1 and P2 by Lemma 4.10. If no arc or a biarc is to be added, then β itself will be an isomorphism from CP1 to CP2 of the required type by the assumption on β.

Now let CP1 and CP2 be obtained by adding a single arc to P1 and P2 respectively.

We have $dt_{P1}(u) = dt_{P2}(β(u))$   (as β is an isomorphism from P1 to P2 as digraphs also).

Also, label of $u$ in P1 = label of the external vertex β($u$) in β(P1) (= P2) (as β preserves labels).

Hence in the completions, $dt_{CP1}(u) = dt_{CP2}(β(u))$ by Observation 4.4(a).

Hence   $dt_{CP1}(u) - dt_{P1}(u) = dt_{CP2}(β(u)) - dt_{P2}(β(u))$.

Thus the change in $dt_{P1}(u)$ due to the new arc added to P1 while completing it is *same as* the change in $dt_{P2}(β(u))$ due to the new arc added to β(P1) = P2 while completing it. Hence the "direction at $u$ of the single unpaired arc" with which $u$ and $w$ are joined while completing P1 is same as the "direction at β(u) of the single unpaired arc" with which β(u) and β(w) are joined while completing P2. Hence β will be an isomorphism of the required type from CP1 to CP2. □

**Corollary 5.3.** A dapasting of a pair of dacards of a graph/digraph G as members of *Dadeck*(G) has a *unique completion*.

**Corollary 5.4.** If a given dapasting of a pair of dacards of a graph/digraph G as members of *Dadeck(G) is* a dapasting in two da-hypomorphs $H_1$ and $H_2$ of G, then $H_1 \cong H_2$.

**Note 5.5.** The analogues of Corollary 5.3 and Corollary 5.4 for pastings of a pair of cards of a digraph are *not true*. In Example 1.4 above, P is a *pasting* of the two cards D−u and D−v of D as members of *Deck*(D) in D as well as in E so that both D and E are completions of P, but D $\not\cong$ E.

**Note 5.6.** *The completions of two nonisomorphic dapastings of two dacards of a graph/digraph can be isomorphic.* In Figure 3 (Example 4.9), P1 and P2 are two nonisomorphic dapastings of a pair of dacards identical to G−a and G−b and both of them occur as dapastings in G; but their completions (both G) are isomorphic. *This does not violate Theorem 5.2* as there is no automorphism of G taking the set of external vertices of the cards in $P_1$ to that in $P_2$. Again, in the case of digraphs $C_p$, $p = 2^n + 2^n$, $n \geq 2$ defined in [11, Theorem 4], the dacards $C_p - v_1$ and $C_p - v_{1+p/2}$ are identical as there is an automorphism of $C_p$ taking $v_1$ to $v_{1+p/2}$. The dapastings of the identical pairs of dacards ($C_p - v_{1+p/4}$, $C_p - v_1$) and ($C_p - v_{1+p/4}$, $C_p - v_{1+p/2}$) as in $C_p$ are not isomorphic as the numbers of arcs in them are different. However, both these dapastings have $C_p$ as a completion by Definition 4.1. The adjacency matrix of digraph $C_8$ is in Figure 5.

```
0 1 1 0 1 0 1 0
0 0 1 1 0 1 0 1
0 0 0 1 1 0 1 0
1 0 0 0 0 1 0 1
1 0 1 0 0 1 1 0
0 1 0 1 0 0 1 1
1 0 1 0 0 0 0 1
0 1 0 1 1 0 0 0
```

**Figure 5:** Adjacency matrix of digraph $C_8$ in which the identical pairs of dacards ($C_8 - v_3$, $C_8 - v_1$) and ($C_8 - v_3$, $C_8 - v_5$) are dapasted nonisomorphically.

The following very crucial definition helps us in giving a common *necessary and sufficient condition* for two da-hypomorphs to be isomorphic in terms of dapastings of cards in them.

**Definition 5.7.** Let G and H be two da-hypomorphic digraphs and let A and B be dacards of G obtained by deleting distinct vertices of G. The dacards A and B are said to be *dapasted isomorphically* in G and H if there are dapastings P1 and P2 of them as members of *Dadeck*(G) in G and in H respectively such that P1 and P2 are isomorphic as dapastings. (There may be dapastings P3 and P4 of A and B in G and in H respectively such that P3 and P4 are not isomorphic.) The corresponding definition holds for pastings of two cards A and B of G when G and H are hypomorphic digraphs/graphs.

**Theorem 5.8.** Let H be a da-hypomorph of a graph/digraph G. Now G ≅ H *if and only if* G has *a pair of dacards* A and B which is dapasted isomorphically in G and H.

**Proof.** *"If" Part*: Dacards A and B of G are dapasted isomorphically in G and H by hypothesis. Hence by the above definition, there exist a dapasting P of A and B *in* G and a dapasting P* of A and B *in* H such that P ≅ P* as dapastings. Hence the completions of P and P*, which are G and H respectively (by Definition 4.2) are isomorphic by Theorem 5.2.

*"Only if" part:* Let G ≅ H and β be an isomorphism from G to H. Let A and B be a pair of dacards in the common dadeck of G and H and let P be a dapasting of them in G with $u$ and $v$ as the external vertices of A and B respectively. Hence A (≅ G−$u$) and B (≅ G−$v$) are isomorphic as dacards to β(G) − β($u$) and β(G) − β($v$) of H respectively with $dt_G(u) = dt_H(β(u))$ and $dt_G(v) = dt_H(β(v))$ and their dapasting P in G (as given in the remark following Definition 4.1), is isomorphic to the dapasting of β(G) − β($u$) (≅ A) and β(G) − β($v$) (≅ B) in H under β. □

**Corollary 5.9.** Let H1 and H2 be a pair of da-hypomorphs of a graph/digraph. A pair of dacards is dapasted isomorphically in H1 and H2 *if and only if* every pair of dacards is pasted isomorphically in H1 and H2.

**Proof.** *"If"* part is obvious.

*"Only if" part:* A pair of dacards is dapasted isomorphically in H1 and H2. Hence H1 ≅ H2 by the theorem. Now let A and B be a pair of dacards in the common dadeck of H1 and H2. As in the proof for the *"only if"* part of the theorem, A and B are dapasted isomorphically in H1 and H2. □

**Note 5.10.** Results analogous to the above two are not true for hypomorphic digraphs. In Figure 1, D and E are hypomorphic and cards D−$u$ and D−$v$ are pasted isomorphically in both, but D ≇ E.

**Theorem 5.11. (a).** Whenever G and H are da-hypomorphs of a graph/digraph, either (1) every pair of dacards in their common dadeck is dapasted isomorphically in G and H, or (2) no pair of dacards in their common dadeck is dapasted isomorphically in G and H.

**(b)**. When (1) holds in Part (a) above, G and H are isomorphic and when (2) holds, G and H are not isomorphic.

**Proof.** Part (a): *Suppose* (2) *is not true*. Now there is a pair of dacards A and B in the common dadeck of G and H such that A and B are dapasted isomorphically in G and H. Hence by Corollary 5.9, every pair of dacards in the common dadeck of G and H is dapasted isomorphically in G and H and (1) holds.

Part (b) follows from Theorem 5.8. □

## 6. Da-reconstructibility of a new class of digraphs

A general theorem on pastings and the da-reconstructibility of a new class of digraphs are proved here.

**Definition 6.1.** Two dacards A and B of a graph/digraph G are said to have a *unique dapasting as members of Dadeck(G)* if in the terminology of Observation 4.14(a), the following are independent of the choice of a dapasting of A and B as members of *Dadeck*(G): (i) the vertices $v^* \in V(A)$ and $u^* \in V(B)$ such that $\beta_B^{-1}\beta_A$ is an isomorphism from $A - v^*$ to $B - u^*$ and (ii) $\beta_B^{-1}\beta_A(w)$ for $w \in V(A) - \{v^*\}$. Unique pasting of two cards of a graph/digraph is defined similarly.

**Definition 6.2.** For a digraph D, the graph having the same vertex set as D in which there is an edge between two vertices if and only if there is at least one arc between them in D is called the *underlying graph* of D and is denoted as U(D). For a given graph G, a digraph D having the same vertex set as G in which there is an arc (an unpaired arc or a biarc) between two vertices if and only if there is an edge between them in G is called a *biorientation* of G.

**Theorem 6.3.** Let G be a graph having a pair $v$, $w$ of distinct vertices such that $A \cong G - v$ and $B \cong G - w$ and *Deck*(G) has no card outside the multiset $\{A, B\}$ which is isomorphic to a member of $\{A, B\}$. Let O(A) and O(B) be the dacards into which A and B transform when G is bioriented as D. If A and B have a *unique* pasting as members of *Deck*(G), then all the dapastings of O(A) and O(B) as members of *Dadeck*(D) are isomorphic.

**Proof.** Cards A and B have a unique pasting P as members of *Deck*(G). Hence by Definition 6.1 there exist unique vertices $v^* \in V(A)$ and $u^* \in V(B)$ and a unique isomorphism $\psi$ from $A - v^*$ to $B - u^*$ such that the isomorphism corresponding to $\beta_B^{-1}\beta_A$ of Observation 4.14(a) in every pasting of A and B as members of *Deck*(D) is $\psi$ itself. (1)

Let O(A) and O(B) denote the dacards of D into which A and B have transformed when G gets bioriented as D. Also *no* dacard of *Dadeck*(D) outside the multiset S = {O(A), O(B)} is isomorphic to a member of S by the hypothesis on A and B. (2)

Let P* be a dapasting of the dacards O(A) and O(B) in D. *If possible let there be a dapasting Q\* of* O(A) *and* O(B) *as members of Dadeck*(D) *which is not isomorphic to P\*.* (3)

Edges of A and B get direction and the digraphs O(A) and O(B) are formed when A and B are pasted as P in G. Hence in Q*, which is not isomorphic to P*, $n–2$ vertices of O(A) must coincide with $n–2$ vertices of O(B) in a way *different* from that in P* by (2). Hence by (3) and Definition 6.1 the isomorphisms corresponding to $\beta_B^{-1}\beta_A$ of Observation 4.14(a) for P* and Q* must be different. Hence if Q denotes the pasting of the cards A and B (which are the underlying graphs of O(A) and O(B) respectively) as members of *Dadeck*(G) obtained from the dapasting Q* by dropping the direction of all arcs, then the isomorphisms corresponding to $\beta_B^{-1}\beta_A$ of Observation 4.14(a) for dapastings P and Q must be different and this contradicts (1). Hence (3) is impossible and the theorem holds. □

**Theorem 6.4.** A digraph D is da-reconstructible if its underlying graph U(D) has a pair of cards (A, B) having a unique pasting as members of *Deck*(U(D)) and *Deck*(U(D)) has no card outside the multiset {A, B} which is isomorphic to a member of {A, B}.

**Proof.** Obviously D is *a biorientation* of U(D). Let O(A) and O(B) denote the dacards of D into which the cards A and B of U(D) have transformed during the biorientation of U(D) into D. Now applying the above theorem to graph U(D) and its biorientation D, all the dapastings of the dacards O(A) and O(B) of D as members of *Dadeck*(D) are isomorphic. Hence dacards O(A) and O(B) are dapasted isomorphically in all da-hypomorphs of D and so all da-hypomorphs of D are isomorphic by Theorem 5.8. □

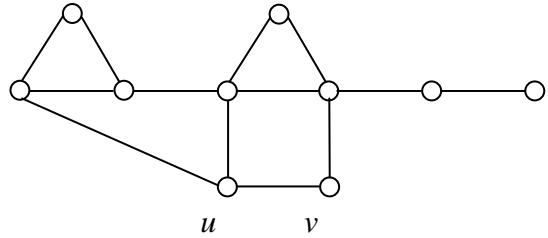

*u*    *v*

**Figure 6**. G (≅ U(D))

**Example 6.5.** All digraphs D with U(D) ≅ G, where G is the graph in Figure 6 are da-reconstructible as below. The graph G (= U(D)) has only one induced subgraph isomorphic to G – *u* – *v*, which has only the identity automorphism. Hence the cards A (≅ G – *u*) and B (≅ G – *v*) in *Deck*(G) = *Deck*(U(D)) satisfy the conditions of Theorem 6.4 and so D is da-reconstructible.

## 7. Known pairs of non-reconstructible digraphs and pastings

For each X ∈ {A, B, C, D, E, F} and each $p = 2^m + 2^n$ with $0 \leq n \leq m$, and $m \geq 1$, Stockmeyer has constructed ([11]) a pair $(X_p, X^*_p)$ of hypomorphic but nonisomorphic digraphs on *p* vertices and these are studied in [10, 12, 13] also. As reported in [13, p. 197], the nine infinite families of pairs in [11], the infinite family of pairs due to Stockmeyer (1988) described as $(G_p, G^*_p)$ in [10], the four pairs (D, E) in Figures 4 and 7 and the three pairs (D, E) in Figures 1 and 8 were the only pairs of non-reconstructible digraphs on ≥ 5 vertices known in 1988. No new pair seems to have been discovered after that ([10, Section 5]).

```
0 1 1 0 1           0 0 0 1 0
0 0 0 1 0           1 0 1 0 1
0 1 0 0 1           1 0 0 0 0
1 0 0 0 0           0 1 0 0 1
0 1 0 1 0           1 0 1 0 0

0 1 1 0 1           0 0 0 1 0
0 0 0 1 0           1 0 1 0 1
0 1 0 1 1           1 0 0 1 0
1 0 1 0 0           0 1 1 0 1
0 1 0 1 0           1 0 1 0 0

0 1 1 1 1 0         0 0 0 0 1 0
0 0 1 1 0 1         1 0 0 0 0 1
0 0 0 1 1 1         1 1 0 0 0 0
0 0 0 0 1 0         1 1 1 0 1 0
0 1 0 0 0 1         0 1 1 0 0 1
1 0 0 1 0 0         1 0 1 1 0 0
      D                   E
```

**Figure 7.** Adjacency matrices of three pairs of hypomorphic but nonisomorphic digraphs

```
      0 1 1 1 1 0              0 0 0 0 0 1
      0 0 1 1 0 1              1 0 0 0 1 0
      0 0 0 1 1 1              1 1 0 0 0 0
      0 0 0 0 0 1              1 1 1 0 1 0
      0 1 0 1 0 0              1 0 1 0 0 1
      1 0 0 0 1 0              0 1 1 1 0 0

0 1 0 1 1 0 0 1          0 0 1 0 0 1 1 0
0 0 1 0 1 1 1 0          1 0 0 1 0 0 0 1
1 0 0 1 0 1 0 1          0 1 0 0 1 0 1 0
0 1 0 0 1 0 1 0          1 0 1 0 0 1 0 1
0 0 1 0 0 1 0 1          1 1 0 1 0 0 1 0
1 0 0 1 0 0 1 0          0 1 1 0 1 0 0 1
1 0 1 0 1 0 0 1          0 1 0 1 0 1 0 0
0 1 0 1 0 1 0 0          1 0 1 0 1 0 1 0
       D                        E
```
**Figure 8.** Adjacency matrices of two pairs of hypomorphic but nonisomorphic digraphs.

**Definition 7.1** ([10], Section 1). Two digraphs D and E are called S-isomorphic (in honor of P.K.Stockmeyer) if either they are isomorphic or there exists a digraph F with a pair of vertices $u$ and $v$ such that there is no arc joining them in F and $D \cong F + uv$ and $E \cong F + vu$.

**Theorem 7.2** ([10], Theorem 4.19). In each of the pairs of hypomorphic but nonisomorphic digraphs (D, E) on five or more vertices known in 1988 other than the four in Figures 4 and 7, the two digraphs in the pair are S-isomorphic.

We prove the following on known pairs of non-reconstructible digraphs.

**Theorem 7.3. (a)** In all but four (those in Figures 4 and 7) of the known pairs of *non-reconstructible digraphs* on $\geq 5$ vertices, the common deck of the digraphs has a pair of cards which is pasted isomorphically in both the digraphs in the pair. **(b)** The four pairs exempted above do not have this property.

**Proof.** *(a) Part*: Let (D, E) be a known pair of hypomorphic but nonisomorphic digraphs other than the four in Figures 4 and 7. By Theorem 7.2 and Definition 7.1, there exists a digraph F with a pair of vertices $u$ and $v$ such that there is no arc joining them in F, $D \cong F + uv$ and $E \cong F + vu$. Hence by Definition 1.3, F with each of $u$ and $v$ relabeled as "$e$" is a pasting of the cards $D - u$ and $D - v$ in D and is a pasting of the cards $E - u$ and $E - v$ in E. Again $D - u \cong F - u \cong E - u$ and $D - v \cong F - v \cong E - v$ and D and E have the same deck. Hence, the cards $F - u$ and $F - v$ in the common deck of D and E are pasted isomorphically (as F with $u$ and $v$ relabeled as "$e$") in D and E.

*(b) part*: This part of the theorem (for the four pairs in Figures 4 and 7) can be directly verified. □

**Observation 7.4.** Only four of the hypomorphic but nonisomorphic pairs of digraphs now known have structure resembling that of a counterexample pair for RC/DARC given by Theorem 5.11.

**Observation 7.5.** None of the existing pairs of hypomorphic but nonisomorphic digraphs on $\geq 5$ vertices is a da-hypomorphic pair.

**Verification.** The digraphs in all the ten known infinite families of non-reconstructible digraphs belong to da-reconstructible families ([10, p. 560]). Hence whenever $(X_p, X^*_p)$ is a pair in one of these ten families, $X_p$ and $X^*_p$ have different dadeck (as $X_p \not\cong X^*_p$ and $X_p$ and $X^*_p$ are da-reconstructible) and so $(X_p, X^*_p)$ is not a da-hypomorphic pair.

The parameter "multiset $\{(dt(v),$ second degree of $w)|\ vw$ is an unpaired arc of D$\}$ of ordered pairs" is a da-reconstructible parameter of a digraph D by Observation 3.7. Among the seven known sporadic pairs (D, E) of non-reconstructible digraphs (those in Figures 1, 4, 7 and 8 above), in all pairs other than the second one in Figure 8, D and E have *different values* for the above parameter. For the second pair in Figure 8, no dacard of D is identical to "the dacard of E obtained by deleting the fourth vertex". (Verified by comparing the ordered pair "(degree triple associated with a dacard, the collection of NDT's (Definition 3.4) of the vertices of that dacard)" corresponding to two different dacards.) Hence in *none* of these seven sporadic pairs, the two digraphs in the pair have the *same dadeck*. □

The above discussion gives rise to the following problem, whose affirmative answer will prove that one type of counterexamples to DARC cannot exist.

**Problem 7.6:** Are all pairs of hypomorphic but nonisomorphic digraphs with their common deck having a pair of cards which is pasted isomorphically in both its member digraphs non-da-hypomorphic?

**8. Conclusion**

1. The claim of the Graph Reconstruction Conjecture (RC) is in terms of a pair of graphs having the same deck. *Pasting of two cards is the key concept* defined in this paper to study such a pair of graphs and it *simplifies* RC and Degree Associated RC as below.

**RC (DARC):** If graphs (digraphs) G and H have the same deck (dadeck), then a pair of cards (dacards) in their common deck (dadeck) is pasted (dapasted) isomorphically in both G and H.

2. Properties of *pastings of two cards* of a graph in a pair of hypomorphs of it and properties of *dapastings of two dacards* of a digraph in a pair of da-hypomorphs of it are alike, whereas properties of *pastings of two cards* of a digraph in a pair of hypomorphs of it are different. (Example: Corollary 5.4 and Note 5.5.) This supports the belief that degree associated reconstruction is the digraph analogue of graph reconstruction.

3. *Optimism*: When G and H are hypomorphs of a graph, they are isomorphic *if and only if a pair of cards* in their common deck *is pasted isomorphically* (Definition 5.7) *in both of them* by Theorem 5.8. This is *not true* for hypomorphs of a digraph (Note 5.10). Stockmeyer calls ([12]) such results *significant* and states that "any proof of the Reconstruction Conjecture must depend on one or more "significant" results in an essential way, since otherwise the proof would hold for digraphs as well". Other results proved here prepare the ground for such an attempt using isomorphic pastings.

4. *Application*: Falsity of RC (DARC) demands the existence of a pair of non-isomorphic graphs (digraphs) with the same deck (dadeck). In such a pair (G, H), no pair of cards (dacards) is pasted (dapasted) isomorphically in both G and H by Theorem 5.11. If we are able to show that a pair of graphs (digraphs) with the above structure cannot exist, then RC (DARC) is proved. One way to do it is to show that a pair of hypomorphic but nonisomorphic digraphs is always non-da-hypomorphic and Problem 7.6 raised above is a part of it. Known pairs of hypomorphic but non-isomorphic digraphs and their properties can serve as scaffolding in such an attempt. Thus the concept *pasting provides a new direction of research* in graph/digraph reconstruction.

**Acknowledgement:** We are thankful to the experts who have given very valuable comments on earlier versions of this manuscript.

**Funding:** Research supported by **CRG/2018/001007/MS** of **SERB (DST)**, **Govt. of India**.

*****